\documentclass[11pt]{amsart}
\usepackage[
paper=a4paper,
text={147mm,230mm},centering
]{geometry}

\iftrue
\makeatletter
\def\@settitle{%
  \vspace*{-20pt}
  \begin{flushleft}%
    \baselineskip14\p@\relax
    \normalfont\bfseries\LARGE
    \@title
  \end{flushleft}%
}
\def\@setauthors{%
  \begingroup
  \def\thanks{\protect\thanks@warning}%
  \trivlist
  \large \@topsep30\p@\relax
  \advance\@topsep by -\baselineskip
  \item\relax
  \author@andify\authors
  \def\\{\protect\linebreak}%
  \authors
  \ifx\@empty\contribs
  \else
    ,\penalty-3 \space \@setcontribs
    \@closetoccontribs
  \fi
  \normalfont
  \endtrivlist
  \endgroup
}
\def\@setabstracta{%
    \ifvoid\abstractbox
  \else
    \skip@25\p@ \advance\skip@-\lastskip
    \advance\skip@-\baselineskip \vskip\skip@
    \box\abstractbox
    \prevdepth\z@ 
    \vskip-10pt
  \fi
}
\renewenvironment{abstract}{%
  \ifx\maketitle\relax
    \ClassWarning{\@classname}{Abstract should precede
      \protect\maketitle\space in AMS document classes; reported}%
  \fi
  \global\setbox\abstractbox=\vtop \bgroup
    \normalfont\small
    \list{}{\labelwidth\z@
      \leftmargin0pc \rightmargin\leftmargin
      \listparindent\normalparindent \itemindent\z@
      \parsep\z@ \@plus\p@
      
    }%
    \item[\hskip\labelsep\bfseries\abstractname.]%
}{%
  \endlist\egroup
  \ifx\@setabstract\relax \@setabstracta \fi
}
\def\section{\@startsection{section}{1}%
  \z@{-1.2\linespacing\@plus-.5\linespacing}{.8\linespacing}%
  {\normalfont\bfseries\large}}
\def\subsection{\@startsection{subsection}{2}%
  \z@{-.8\linespacing\@plus-.3\linespacing}{.3\linespacing\@plus.2\linespacing}%
  {\normalfont\bfseries}}
\def\subsubsection{\@startsection{subsubsection}{3}%
  \z@{.7\linespacing\@plus.1\linespacing}{-1.5ex}%
  {\normalfont\itshape}}
\def\@secnumfont{\bfseries}
\makeatother
\fi 

\usepackage{amssymb}
\usepackage{mathrsfs}
\usepackage[all]{xy}
\usepackage{graphicx,color,float}
\usepackage[bookmarks]{hyperref}
\usepackage{pinlabel}
\usepackage[textwidth=1.3in,color=yellow]{todonotes}
\usepackage{amsmath}
\usepackage{enumitem}
\usepackage{pb-diagram,pb-xy}
\usetikzlibrary{calc}
\usepackage{youngtab}
\usepackage[boxsize=.3 em]{ytableau}

\textwidth=\paperwidth \advance\textwidth by-2\oddsidemargin
\advance\oddsidemargin by-1in
\evensidemargin=\oddsidemargin
\calclayout

\theoremstyle{plain}
\newtheorem{theorem}{Theorem}

\newtheorem{proposition}[theorem]{Proposition}
\newtheorem{lemma}[theorem]{Lemma}
\newtheorem{corollary}[theorem]{Corollary}

\theoremstyle{definition}

\theoremstyle{remark}
\newtheorem{remark}[theorem]{Remark}

\newcommand{\sslash}{\mathbin{/\mkern-6mu/}}

\newcommand{\C}{\mathbb{C}}

\newcommand{\R}{\mathbb{R}}

\newcommand{\Z}{\mathbb{Z}}

\newcommand{\CP}{\mathbb{C}P}

\def\mcal{\mathcal}
\def\frak{\mathfrak}
\def\scr{\mathscr}

\def\to{\mathchoice{\longrightarrow}{\rightarrow}{\rightarrow}{\rightarrow}}
\makeatletter
\newcommand{\shortxra}[2][]{\ext@arrow 0359\rightarrowfill@{#1}{#2}}
\def\longrightarrowfill@{\arrowfill@\relbar\relbar\longrightarrow}
\newcommand{\longxra}[2][]{\ext@arrow 0359\longrightarrowfill@{#1}{#2}}

\makeatother

\begin{document}                                                                          
\title[Chekanov torus and Gelfand--Zeitlin torus]{Chekanov torus and Gelfand--Zeitlin torus in $S^2 \times S^2$}

\author{Yoosik Kim}
\address{Department of Mathematics, Pusan National University}
\email{yoosik@pusan.ac.kr}


\begin{abstract}
The Chekanov torus was the first known \emph{exotic} torus, a monotone Lagrangian torus that is not Hamiltonian isotopic to the standard monotone Lagrangian torus. We explore the relationship between the Chekanov torus in $S^2 \times S^2$ and a monotone Lagrangian torus that had been introduced before Chekanov's construction \cite{Chekanov}. We prove that the monotone Lagrangian torus fiber in a certain Gelfand--Zeitlin system is Hamiltonian isotopic to the Chekanov torus in $S^2 \times S^2$.
\end{abstract}

\maketitle
\setcounter{tocdepth}{1} 

\section{Introduction}

A symplectic manifold $(X, \omega)$ is called \emph{monotone} if the class of the symplectic form $\omega$ and the first Chern class of the tangent bundle $TX$ in $H^2(X; \R)$ are positively proportional, i.e. $[\omega] = r \cdot c_1 (TX)$ for a positive real number $r > 0$. Each Lagrangian submanifold $L$ of $X$ carries two group homomorphisms$\colon$
\begin{itemize}
\item $I_\omega \colon \pi_2(X,L) \to \R$ is given by the symplectic area,
\item $I_\mu \colon \pi_2(X,L) \to \Z$ is given by the Maslov index.
\end{itemize}
A Lagrangian submanifold $L$ is said to be \emph{monotone} if there exists a positive real number $s > 0$ such that $I_\omega (\beta) = s \cdot I_\mu (\beta)$ for every $\beta \in \pi_2(X,L)$. This constant $s$ is called the \emph{monotonicity} of $L$. The notion of a monotone Lagrangian submanifold was introduced by Y.-G. Oh to construct its Floer cohomology in \cite{OhMonotone, OhMonotone2}.

As an attempt to classify the monotone Lagrangians in a symplectic manifold, it has been an interesting problem to construct monotone Lagrangian tori that are not related by any Hamiltonian isotopy. In the symplectic vector space $\R^{2n}$, Chekanov \cite{Chekanov} constructed a monotone Lagrangian torus that is \emph{not} Hamiltonian isotopic to any standard product torus. Moreover, by embedding a suitable subset of $\R^{2n}$ symplectically, one can produce  a monotone Lagrangian torus by transporting a Chekanov exotic torus of $\R^{2n}$ into a closed symplectic manifold including the complex projective space or the product of projective lines. Such an embedded monotone Lagrangian torus is also called a \emph{Chekanov torus}. A monotone torus not Hamiltonian isotopic to the standard one is called an \emph{exotic torus}. The Chekanov torus is the first known exotic monotone Lagrangian torus. Recently, there have been exciting developments constructing infinitely many exotic Lagrangian tori, see \cite{Viannacp2infi, AurouxExotic, ViannadelP} for instance.

In this article, we are concerned with the Chekanov torus in $\CP^1 \times \CP^1$. To begin with, we briefly recall several ways of constructing the Chekanov torus in $\CP^1 \times \CP^1 \simeq S^2 \times S^2$.

\begin{itemize}
\item (Chekanov--Schlenk torus $T_{\mathrm{CS}}$ in \cite{ChekanovSchlenk}) Let $\mathbb{D}(\sqrt{2}) := \{ \zeta \in \C \mid | \zeta | < \sqrt{2} \}$ be an open disk in the complex plane with radius $\sqrt{2}$ and let $\mathbb{H}(\sqrt{2}) := \{ \zeta \in \mathbb{D}(\sqrt{2}) \mid \mathrm{Im} (\zeta) > 0\}$ be an open half-disk. Choose any simple closed curve $\Gamma$ in $\mathbb{H}(\sqrt{2})$ which bounds the region having the area $\pi/2$. The product space carries the Hamiltonian $S^1$-action given by
$$
e^{\sqrt{-1} \theta} \cdot (\zeta_1, \zeta_2) \to \left(e^{-\sqrt{-1}\theta} \zeta_1, e^{\sqrt{-1}\theta} \zeta_2 \right).
$$

Consider the product space $\CP^1 \times \CP^1$ of projective planes equipped with the product symplectic form $2 (\omega_{\mathrm{FS}} \oplus \omega_{\mathrm{FS}})$ where $([v_0 : v_1], [w_0 : w_1])$ is its homogeneous coordinate. Recall that $2 \omega_{\mathrm{FS}} = \frac{1}{2} \omega_{\mathrm{std}}$ and the symplectic area of $\CP^1$ measured by $2 \omega_{\mathrm{FS}}$ is $2\pi$. We denote the diagonal map by
\begin{equation}\label{equ_diagmap}
\Delta_{\mathbb{D}} \colon \mathbb{D}(\sqrt{2}) \to \mathbb{D}(\sqrt{2})  \times \mathbb{D}(\sqrt{2}).
\end{equation}
Also, the composition of the inclusion $\mathbb{H}(\sqrt{2}) \to \mathbb{D}(\sqrt{2})$ and the diagonal map $\Delta_{\mathbb{D}}$ is denoted by
\begin{equation}\label{equ_diagmaph}
\Delta_{\mathbb{H}} \colon \mathbb{H}(\sqrt{2}) \to \mathbb{D}(\sqrt{2}) \to \mathbb{D}(\sqrt{2})  \times \mathbb{D}(\sqrt{2}).
\end{equation}
We have a symplectic embedding
\begin{equation}\label{equ_symplecticembeddingrho}
\rho \colon \left( \mathbb{D}(\sqrt{2}) \times \mathbb{D}(\sqrt{2}), ( \omega_{\mathrm{std}} \oplus \omega_{\mathrm{std}}) |_{\mathbb{D}(\sqrt{2}) \times \mathbb{D}(\sqrt{2})} \right) \to \left( \CP^1 \times \CP^1, 2  (\omega_{\mathrm{FS}} \oplus \omega_{\mathrm{FS}}) \right)
\end{equation}
such that the image of $\rho$ is the intersection of $v_0 \neq 0$ and $w_0 \neq 0$. Chekanov and Schlenk in \cite{ChekanovSchlenk} constructed the monotone Lagrangian torus $T_{\mathrm{CS}}$ as follows$\colon$
\begin{equation}\label{equ_TCS}
T_{\mathrm{CS}} := \left\{ \rho \left(e^{-\sqrt{-1}\theta} \zeta, e^{\sqrt{-1}\theta} \zeta \right) \in \CP^1 \times \CP^1 \mid \theta \in [0, 2\pi], \zeta \in \Gamma   \right\}.
\end{equation}
\end{itemize}

\begin{itemize}
\item (Entov--Polterovich torus $T_{\mathrm{EP}}$ in \cite{EntovPolterovich}) Let us regard the product space $S^2 \times S^2$ as an embedded submanifold of $\R^3 \times \R^3$ equipped with $\frac{1}{2}(\omega_{\mathrm{std}} \oplus \omega_{\mathrm{std}})$. 
Entov and Polterovich in \cite{EntovPolterovich} constructed a monotone Lagrangian torus as follows$\colon$
$$
T_{\mathrm{EP}} := \left\{ (\mathbf{a}, \mathbf{b}) \in S^2 \times S^2 \mid (\mathbf{a}+\mathbf{b}) \cdot \mathbf{e}_1 = 0, \mathbf{a} \cdot \mathbf{b} = - 1/2 \right\}
$$
where $\mathbf{e}_1 = (1, 0, 0)$.

\item (Fukaya--Oh--Ohta--Ono torus $T_{\mathrm{FOOO}}$ in \cite{FOOOS2S2}) Start with the symplectic toric orbifold associated to the triangle whose vertices are $(0,0), (0,1)$, and $(2,0)$ as the moment polytope. By replacing a neighborhood of the singular point with the Milnor fiber to obtain a symplectic manifold isomorphic to $S^2 \times S^2$. Through this process, one obtains a semi-toric system $\Phi_{\mathrm{semi}}$ whose image is the above polytope. Fukaya--Oh--Ohta--Ono monotone Lagrangian torus is located at the center $(1/2, 1/2)$, that is,
$$
T_{\mathrm{FOOO}} := \Phi^{-1}_{\mathrm{semi}} (1/2, 1/2).
$$

\item (Albers--Frauenfelder torus $T_{\mathrm{AF}}$ in \cite{AlbersFrauenfelder}) Let $\Delta$ be the diagonal of $S^2 \times S^2$. Choose a symplectomorphism
$$
\upsilon \colon (S^2 \times S^2) \backslash \Delta \to D_1^* S^2
$$
where $D_1^* S^2$ is the open unit disk bundle. Albers--Frauenfelder monotone torus is defined by 
$$
T_{\mathrm{AF}} := \upsilon^{-1} \left( \left\{ (\mathbf{p},\mathbf{q}) \in D_1^* S^2 \mid |\mathbf{p}| = 1/2, (\mathbf{p} \times \mathbf{q}) \cdot \mathbf{e}_1 = 0\right\} \right). 
$$

\item (Biran--Cornea torus $T_{\mathrm{BC}}$ in \cite{BiranCornea})
Let us start with the diagonal $\Delta$ of $\CP^1 \times \CP^1$, a complex hypersurface of $\CP^1 \times \CP^1$. Set 
$$
P_\Delta = \left\{ (\mathbf{a},\mathbf{b}) \in \R^3 \times \R^3 \colon |\mathbf{a}| = |\mathbf{b}| = 1,\, \mathbf{a} \cdot \mathbf{b} = 0 \right\}.
$$
with the circle action 
$$
e^{\sqrt{-1} \theta} \cdot (\mathbf{a},\mathbf{b}) = (\mathbf{a}, (\cos \theta) \mathbf{b} + (\sin \theta) (\mathbf{a} \times \mathbf{b})).
$$
Identifying $\Delta$ with $S^2$, we then have the principal $S^1$-bundle $\pi \colon P_\Delta \to \Delta$ given by $\pi ([(\mathbf{a},\mathbf{b})]) = \mathbf{a}$, whose connection $1$-from for $P_\Delta$ is defined by $\alpha_{(\mathbf{a}, \mathbf{b})} (\mathbf{p}, \mathbf{q}) = \mathbf{q} \cdot (\mathbf{a} \times \mathbf{b})$. According to a work of Biran in \cite{BiranLbse}, there is a symplectomorphism from 
$$
D_{\sqrt{2}}(P_\Delta) = \frac{P_{\Delta} \times D(\sqrt{2})}{(e^{\sqrt{-1}\theta} \cdot (\mathbf{a},\mathbf{b}), \zeta) \sim ((\mathbf{a},\mathbf{b}), e^{\sqrt{-1} \theta} \cdot \zeta)}
$$
adorned with a certain symplectic form defined by using the connection $1$-form to $(S^2 \times S^2) \backslash \overline{\Delta}$ where $\overline{\Delta}$ is the anti-diagonal.

For each monotone Lagrangian submanifold $L$ in $\Delta$ and each $r$ with $0 < r < \sqrt{2}$, we have the submanifold
$$
L_{(r)} = \left\{ \left[ ((\mathbf{a},\mathbf{b}), \zeta ) \right] \in D_{\sqrt{2}}(P_{\Delta}) \mid |\zeta| = r, \pi ([((\mathbf{a},\mathbf{b}), \zeta)]) \in L \right\}.
$$
According to \cite{BiranCornea}, there is a unique value $r$ such that the embedding of $L_{(r)}$ becomes a monotone Lagrangian submanifold in $S^2 \times S^2$. In this case, we choose a great circle of $S^2 \simeq \Delta$ for $L$ and one for $r$. The choice leads to Brian--Cornea monotone Lagrangian torus
$$
T_{\mathrm{BC}} = L_{(1)}
$$
in $S^2 \times S^2$. For an explicit description for the embedding and the embedded monotone Lagrangian, the reader is referred to \cite[Section 4]{OakleyUsher}.  
\end{itemize}

The relation between the above monotone tori had been explored. In \cite{FOOOS2S2}, Fukaya, Oh, Ohta, and Ono discussed the relation between $T_{\mathrm{FOOO}}$ and $T_{\mathrm{AF}}$. Gadbled \cite{Gadbled} showed that the $T_{\mathrm{BC}}$ and $T_{\mathrm{CS}}$ are Hamiltonian isotopic. Oakley and Usher ultimately proved that all five monotone Lagrangian tori are Hamiltonian isotopic to each other by constructing detailed and explicit symplectomorphisms in \cite{OakleyUsher}. 

\begin{theorem}[Theorem 1.1 in \cite{OakleyUsher}]
The five monotone Lagrangian tori $T_{\mathrm{CS}}$, $T_{\mathrm{EP}}$, $T_{\mathrm{FOOO}}$, $T_{\mathrm{AF}}$, and $T_{\mathrm{BC}}$ listed above are Hamiltonian isotopic to each other.
\end{theorem}

The main goal of this article is to add one more monotone Lagrangian torus to the above list. This torus had been constructed even before the Chekanov's construction of ``twist tori" in \cite{Chekanov}. Specifically, we shall prove that the monotone Gelfand--Zeitlin (GZ) Lagrangian torus fiber in the orthogonal Grassmannian $\mathrm{OG}(1, \C^4)$ is Hamiltonian isotopic to one (and hence each) of the above listed tori. 

Let us recall the quickest way of constructing the monotone GZ Lagrangian torus.
\begin{itemize}
\item (Gelfand--Zeitlin torus $T_{\mathrm{GZ}}$ in \cite{GuilleminSternbergCCI}) Consider the set $\mcal{O}$ of $(4 \times 4)$ skew-symmetric matrices with prescribed four eigenvalues $\pm \lambda \sqrt{-1}$ and $0$ with multiplicity two for some real number $\lambda$. In other words, the space $\mcal{O}$ is the orbit of the block diagonal matrix
\begin{equation}\label{equ_blockdiag}
\begin{bmatrix}
0 & \lambda & 0 & 0\\
-\lambda & 0 & 0 & 0\\
0 & 0 & 0 & 0\\
0 & 0 & 0 & 0
\end{bmatrix}
\end{equation}
under the adjoint $\mathrm{SO}(4)$-action. It carries the Kirillov--Kostant--Souriau (KKS) symplectic form $\omega_{\mathrm{KKS}}$, which is determined by $\lambda$. cf. Remark~\ref{remark_kksformconv} below. To match up with the convention on the symplectic form $2 (\omega_{\mathrm{FS}} \oplus \omega_{\mathrm{FS}})$ on $\CP^1 \times \CP^1$, we take $\lambda = 1$.

For a matrix $A \in \mcal{O}$, let $A^{(k)}$ be the leading principal submatrix of the size $(k \times k)$. Since $A$ is skew-symmetric, $A^{(1)}$ is the zero matrix. In order to define the desired GZ torus, we impose the following conditions on $A^{(2)}$ and $A^{(3)}$.
\begin{enumerate}
\item The eigenvalues of $A^{(2)}$ are required to be $0$ with multiplicity two. 
\item The eigenvalues of $A^{(3)}$ are required to be $\pm \frac{1}{2}\sqrt{-1}$, and $0$. 
\end{enumerate}
The conditions $(\mathrm{1})$ and $(\mathrm{2})$ determine the GZ Lagrangian torus $T_{\mathrm{GZ}}$ in $\mcal{O}$. That is,
\begin{equation}\label{equ_TGZdescr}
T_{\mathrm{GZ}} = \left\{ A \in \mcal{O} ~\Big{|}~ A^{(2)} = O \mbox{ and the eigenvalues of $A^{(3)}$ are $\pm \frac{1}{2}\sqrt{-1}, 0$} 
\right\}.
\end{equation}
\end{itemize}

\begin{remark}\label{remark_kksformconv}
When defining the KKS form $\omega_{\mathrm{KKS}}$, the symplectic form is sometimes normalized by dividing the form by $2\pi$ so that the class of the symplectic form agrees with an \emph{integral} multiple of the Chern class. In this article, we do \emph{not} normalize the KKS form to facilitate comparison of the KKS form with other symplectic forms. Namely, it is defined by 
$$
\omega_{\mathrm{KKS}} \left( \mathrm{ad}_{\xi_1} (A), \mathrm{ad}_{\xi_2} (A) \right) = \langle A, [\xi_1, \xi_2] \rangle
$$ 
where $A \in \mcal{O}$ and $\mathrm{ad}_{\xi_i} (A)$ is a tangent vector at $A$ of $\mcal{O}$ for $i = 1, 2$.
\end{remark}

Additional explanation on $T_{\mathrm{GZ}}$ and $\mcal{O}$ is in order. The real special orthogonal group $\mathrm{SO}(4) := \mathrm{SO}(4; \R)$ transitively acts on $\mcal{O}$ by conjugation. A choice of an element of $\mcal{O}$ gives rise to a diffeomorphism
$
\mcal{O} \simeq {\mathrm{SO}(4)}/{\mathrm{S}(\mathrm{O}(2) \times \mathrm{O}(2))}$.
Therefore, the orbit $\mcal{O}$ is diffeomorphic to the orthogonal Grassmannian $\mathrm{OG}(1,\C^4)$, which is the space of isotropic subspaces of $\C^4$ with respect to a non-degenerate symmetric bilinear form. Thus, $\mathrm{OG}(1,\C^4)$ can be regarded as a quadric hypersurface. As any two smooth quadric hypersurfaces are symplectomorphic, the quadric hypersurface is symplectomorphic to the Segre variety $\left\{ [\mathbf{z}] \in \CP^3 \mid z_0 z_1 - z_2 z_3 = 0 \right\}$. Hence, the space $\mcal{O}$ is regarded as the product $\CP^1 \times \CP^1$ of two complex projective planes.

The min-max principle, a system of inequalities between eigenvalues of submatrices of a skew-symmetric matrix, says that the eigenvalues of $A^{(2)}$ and $A^{(3)}$ are of the form$\colon$
\begin{itemize}
\item the eigenvalues of $A^{(2)}$ are $\pm \lambda^{(2)} \sqrt{-1}$.
\item the eigenvalues of $A^{(3)}$ are $\pm \lambda^{(3)} \sqrt{-1}, 0$.
\end{itemize}
Each real-valued function determined by $\lambda^{(\bullet)}$ generates a Hamiltonian circle action (on an open dense subset of $\mcal{O}$) by the result of Guillemin--Sternberg \cite{GuilleminSternbergGC}. Indeed, the fiber $T_{\mathrm{GZ}}$ is a free Hamiltonian $T^2$-orbit so that it is a Lagrangian torus. Employing the technique of the gradient Hamiltonian disks in \cite{ChoKimMONO}, the torus $T_{\mathrm{GZ}}$ was shown to be monotone.

\begin{lemma}[Proposition 3.7 in \cite{KimQuad}]\label{lemma_GZmonotone}
The above Gelfand--Zeitlin torus fiber $T_{\mathrm{GZ}}$ is monotone.
\end{lemma}

The main theorem of this article is stated below.

\begin{theorem}\label{theorem_main}
There is a symplectomorphism from $\mcal{O}$ to $\CP^1 \times \CP^1$ taking the Gelfand--Zeitlin torus $T_\mathrm{GZ}$ to the Chekanov--Schlenk torus $T_\mathrm{CS}$. Consequently, all monotone Lagrangian tori $T_{\mathrm{CS}}$, $T_{\mathrm{EP}}, T_{\mathrm{FOOO}}, T_{\mathrm{AF}}$, $T_{\mathrm{BC}}$, and $T_{\mathrm{GZ}}$ are Hamiltonian isotopic to each other.
\end{theorem}

There had been decisive clues that $T_{\mathrm{GZ}}$ is Hamiltonian isotopic to the Chekanov torus in literatures. First, Nishinou--Nohara--Ueda \cite{NishinouNoharaUeda2} described and constructed the  toric degeneration of the Gelfand--Zeitlin system. The image of the system on the quadric hypersurface of complex dimension two agrees with the moment polytope of $\CP(1,1,2)$, which is one of the starting points of this work. Second, the disk potential of the GZ torus fiber $T_{\mathrm{GZ}}$ computed in \cite{KimQuad} agrees with that of $T_{\mathrm{CS}}$ in \cite{AurouxTdual, ChekanovSchlenk} and that of $T_{\mathrm{FOOO}}$ in \cite{FOOOS2S2} up to some coordinate changes. Theorem~\ref{theorem_main} claims a stronger relationship between $T_{\mathrm{GZ}}$ and $T_\mathrm{CS}$.

\section{Proof of Theorem~\ref{theorem_main}}

Consider the product space $\CP^1 \times \CP^1$ of complex planes equipped with the product symplectic form $2 (\omega_{\mathrm{FS}} \oplus \omega_{\mathrm{FS}})$ where $\omega_{\mathrm{FS}}$ is the Fubini--Study form on $\CP^1$. It has the $S^1$-action given by
\begin{equation}\label{equ_s11}
([v_0 : v_1], [w_0 : w_1]) \mapsto \left( \left[ v_0 : e^{-\sqrt{-1} \theta} v_1 \right], \left[w_0 : e^{\sqrt{-1} \theta} w_1\right] \right).
\end{equation}

The product space $\CP^1 \times \CP^1$ is embedded into $\CP^3$ via the Segre embedding
$$
\sigma \colon \CP^1 \times \CP^1 \to \CP^3 \quad ([v_0 : v_1], [w_0 : w_1]) \mapsto [v_0 w_1 : v_1 w_0: v_0 w_0: v_1 w_1].
$$
The image under the Segre map $\sigma$ is a hypersurface of $\CP^3$ defined by $z_0 z_1 - z_2 z_3 = 0$ where $[\mathbf{z}] := [z_0: z_1: z_2: z_3]$ is the homogeneous coordinates for $\CP^3$. The hypersurface is denoted by
\begin{equation}\label{equ_x0x1x2x3}
\mcal{Q}^\prime := \left\{[\mathbf{z}] \in \CP^3 \mid z_0 z_1 - z_2 z_3 = 0 \right\}.
\end{equation}

Let $S^1$ act on $\CP^3$ as follows. 
\begin{equation}\label{equ_s12}
[z_0: z_1: z_2: z_3] \mapsto \left[e^{\sqrt{-1} \theta} z_0: e^{-\sqrt{-1} \theta} z_1: z_2: z_3\right].
\end{equation}
Adorning $\CP^3$ with twice Fubini--Study form $2 \omega_{\mathrm{FS}}$, the $S^1$-action becomes Hamiltonian. Setting $\| \mathbf{z} \|^2 :=  {|z_0|^2 + |z_1|^2 + |z_2|^2 + |z_3|^2}$, we choose a moment map $\mu_{S^1} \colon \CP^3 \to \R$ with respect to the symplectic form $2 \omega_{\mathrm{FS}}$ for the action~\eqref{equ_s12} as 
\begin{equation}\label{equ_momentmapz}
\mu_{S^1} (\mathbf{z}) =  \frac{|z_0|^2 - |z_1|^2}{\| \mathbf{z} \|^2}.
\end{equation}
Then the map $\sigma$ is $S^1$-equivariant. The embedded variety $\mcal{Q}^\prime$ is equipped with the K{\"a}hler form inherited from $(\CP^3, 2 \omega_{\mathrm{FS}})$. We denote the form restricted to $\mcal{Q}^\prime \subset \CP^3$ by $\omega_{\mcal{Q}^\prime}$.

To show that the map $\sigma$ preserves the symplectic form, it suffices to show that it is a symplectomorphism on an affine chart, an open dense subset of $\CP^1 \times \CP^1$. A straightforward computation shows that $\sigma$ preserves the symplectic form.

\begin{proposition}\label{Proposition_Segreembedding}
The Segre embedding $\sigma$ is an $S^1$-equivariant symplectomorphism.
\end{proposition}

Under the coordinate change
\begin{equation}\label{equ_coordinatechange}
\begin{bmatrix}
z_0 \\ z_1 \\ z_2 \\ z_3
\end{bmatrix}
=
\begin{bmatrix}
1 & \sqrt{-1} & 0 & 0\\
1 & -\sqrt{-1} & 0 & 0\\
0 & 0 & 1 & \sqrt{-1}\\
0 & 0 & -1 & \sqrt{-1}
\end{bmatrix}
\begin{bmatrix}
x_0 \\ x_1 \\ x_2 \\ x_3
\end{bmatrix},
\end{equation}
the variety $\mcal{Q}^\prime$ in~\eqref{equ_x0x1x2x3} maps into the Fermat hypersurface$\colon$
\begin{equation}\label{equ_q2}
\mcal{Q} = \left\{ [x_0 : x_1 : x_2 : x_3] \in \CP^3 \mid x_0^2 + x_1^2 + x_2^2 + x_3^2 = 0 \right\}.
\end{equation}
The linear transformation~\eqref{equ_coordinatechange} is denoted by $\Lambda \colon \mcal{Q} \to \mcal{Q}^\prime$. 

Let $G := \mathrm{SO}(2) = \mathrm{SO}(2;\R) \simeq \mathrm{diag}\left(\mathrm{SO}(2;\R), I_{2} \right)$ act on $\C^4$ linearly. The linear $\mathrm{SO}(2)$-action induces the action on $\CP^3$ and the action on $\mcal{Q}$. Let us fix the following isomorphism
\begin{equation}\label{equation_SO2S1}
\iota \colon S^1 \to \mathrm{SO}(2; \R), \quad e^{\sqrt{-1} \theta} \mapsto 
\begin{bmatrix}
\cos \theta & - \sin \theta \\
\sin \theta & \cos \theta 
\end{bmatrix}.
\end{equation}
The hypersurface $\mcal{Q}$ in~\eqref{equ_q2} admits the $S^1$-action via~\eqref{equation_SO2S1}. We then observe the following.

\begin{lemma}\label{lemma_qq'eq}
The map $\Lambda \colon \mcal{Q} \to \mcal{Q}^\prime$ is $S^1$-equivariant.
\end{lemma}

\begin{proof}
Using the identification~\eqref{equation_SO2S1} and the expression~\eqref{equ_coordinatechange}, one can directly verify the $S^1$-equivalence of $\Lambda$.
\end{proof}

We now discuss symplectic forms on $\mcal{Q}$. Two adorned symplectic forms on $\mcal{Q}$ are taken into account. The first one is the pull-backed symplectic form $\Lambda^* \omega_{\mcal{Q}^\prime}$. The second one is the symplectic form coming from the twice Fubini--Study form $2 \omega_{\mathrm{FS}}$ on the ambient space $\CP^3$. Namely, we take the reduction $\CP^3 \simeq (\C^4 - \{0\}) \sslash S^1$ where the $S^1$-action is diagonal action with the stability condition $\| \mathbf{x} \| = \sqrt{2}$ so that $\CP^3$ is equipped with twice Fubini--Study form $2 \omega_{\mathrm{FS}}$. The K{\"a}hler form restricted to $\mcal{Q}$ is denoted by $\omega_\mcal{Q}$. 

The Fermat hypersurface $\mcal{Q}$ carries two circle actions.
First, the moment map of the corresponding $S^1$-action on $\mcal{Q}$ with respect to the first symplectic form $\Lambda^* \omega_{\mcal{Q}^\prime}$ is obtained by replacing $\mathbf{z}$ with $\Lambda(\mathbf{x})$ in~\eqref{equ_momentmapz} as follows.
\begin{equation}\label{equ_momentq2psiwz}
\mu_{S^1}(\mathbf{x}) =   \frac{\sqrt{-1} \left(\overline{x}_0 x_1 - x_0\overline{x}_1 \right)}{\| \mathbf{x} \|^2}.
\end{equation}
Second, by taking the Killing form
$$
\langle \xi_1, \xi_2 \rangle = -\frac{1}{2} \mathrm{Tr}(\xi_1, \xi_2),
$$
we identify the Lie algebra $\frak{g}$ and its dual Lie algebra $\frak{g}^*$. Under this identification, a moment map of the above $\mathrm{SO}(2)$-action on $\mcal{Q}$ is then of the form
\begin{equation}\label{equ_Gactionqq}
\mu_G \colon \mcal{Q} \to \frak{g}, \quad \mu_G(\mathbf{x}) = \frac{\sqrt{-1}}{ \| \mathbf{x} \|^2 } \begin{bmatrix}
0 & (\overline{x}_0 x_1 - x_0\overline{x}_1) \\
- (\overline{x}_0 x_1 - x_0\overline{x}_1) & 0 
\end{bmatrix}.
\end{equation}
Via the isomorphism~\eqref{equation_SO2S1}, the $\mathrm{SO}(2)$-action can be regarded as the $S^1$-action. By the functoriality of moment maps, a moment map of the $S^1$-action can be chosen as
\begin{equation}\label{equ_momentpsix}
(\iota^* \circ \mu_{G})(\mathbf{x}) =  \frac{\sqrt{-1} \left( \overline{x}_0 x_1 - x_0\overline{x}_1 \right)}{\| \mathbf{x} \|^2} \colon \mcal{Q} \to \R.
\end{equation}

The following lemma compares those two circle actions.

\begin{lemma}
The Hamiltonian function~\eqref{equ_momentq2psiwz} on $\left( \mcal{Q}, \Lambda^* \omega_{\mcal{Q}^\prime} \right)$ and the Hamiltonian function~\eqref{equ_momentpsix} on $\left( \mcal{Q}, \omega_\mcal{Q} \right)$ generate the same $S^1$-action on $\mcal{Q}$.
\end{lemma}

\begin{proof}
We shall compare those two $S^1$-actions on $\mcal{Q}$ after passing them to $\mcal{Q}^\prime$ via the isomorphism $\Lambda$ in~\eqref{equ_coordinatechange}. Recall that the $S^1$-action on $\mcal{Q}^\prime$ corresponding to the first $S^1$-action is described in~\eqref{equ_s12}. By Lemma~\ref{lemma_qq'eq}, the $S^1$-action on $\mcal{Q}^\prime$ corresponding to the second $S^1$-action is exactly the action in~\eqref{equ_s12} as desired.
\end{proof}

To interpolate these forms, we need an equivariant version of the Moser theorem, which is stated below.

\begin{theorem}[Equivariant Moser Theorem]\label{thm_equimoserthm}
Let $X$ be a compact symplectic manifold. Suppose that $\omega_0$ and $\omega_1$ are two symplectic forms in the same cohomology class, that is, $[\omega_0] = [\omega_1]$ in $H^2(X; \C)$. If $\omega_t := (1-t) \omega_0 + t \omega_1$ is symplectic for each $t \in [0,1]$, then there exists an isotopy
$$
\phi \colon [0,1] \times X \to X
$$
such that $\phi_t^* \omega_t = \omega_0$ for each $t$. 

If a compact Lie group $G$ acts on $X$ symplectically with respect to $\omega_t$ for each $t$, then the map $\phi_t$ is $G$-equivariant for each $t$.
\end{theorem}

Consider $\mcal{Q}$ together with two symplectic forms $\Lambda^* \omega_{\mcal{Q}^\prime}$ and $\omega_{\mcal{Q}}$. Since two forms are K\"{a}hler forms, each form $\omega_t := (1-t)\Lambda^* \omega_{\mcal{Q}^\prime} + t \omega_{\mcal{Q}}$ in the linear interpolation is also symplectic. Since the $S^1$-action is Hamiltonian (and hence symplectic) with respect to both $\omega_0$ and $\omega_1$, the $S^1$-action is symplectic with respect to $\omega_t$ for each $t \in [0,1]$. By Theorem~\ref{thm_equimoserthm}, there exists an $S^1$-equivariant isotopy $\phi_t \colon \mcal{Q} \to \mcal{Q}$ such that $\phi_t^*(\omega_t) = \omega_0$. Set
$$
\Phi := \phi_1 \circ \Lambda^{-1} \circ \sigma  \colon ( \CP^1 \times \CP^1, 2(\omega_{\mathrm{FS}} \oplus \omega_{\mathrm{FS}}))
 \to (\mcal{Q}^\prime, \omega_{\mcal{Q}^\prime} ) \to (\mcal{Q}, \Lambda^* \omega_{\mcal{Q}^\prime}) \to (\mcal{Q},  \omega_{\mcal{Q}}).
$$

\begin{lemma}\label{lemmaxi1equ}
The map $\Phi$ is an $S^1$-equivariant symplectomorphism
\end{lemma}

\begin{proof}
By Proposition~\ref{Proposition_Segreembedding}, the Segre map $\sigma$ is an $S^1$-equivariant symplectomorphism. By Lemma~\ref{lemma_qq'eq}, so is the second map $\Lambda^{-1}$. It follows from Theorem~\ref{thm_equimoserthm} that $\phi_1$ is an $S^1$-equivariant symplectomorphism. Therefore, Lemma~\ref{lemmaxi1equ} is established.
\end{proof}

Note that the quadric $\mcal{Q}$ is the zero locus of the symmetric bilinear form $\mcal{B}$ corresponding to the identity matrix $I_4$ on $\CP^3$. On the other hand, the vanishing condition imposes the isotropy condition on the bilinear form $\mcal{B}$ on $\C^4= \C \langle x_1, x_2, x_3, x_4 \rangle$. Namely, $\mcal{Q}$ is isomorphic to the orthogonal Grassmannian $\mathrm{OG}(1, \C^4)$, which consists of one dimensional subspaces $V$ of $\C^4$ satisfying $\mcal{B}(v_1,v_2)= 0$ for all $v_1, v_2 \in V$. The orthogonal Grassmannian $\mathrm{OG}(1, \C^4)$ is adorned with the complex structure from the identification with $\mathrm{SO}(4; \C)/{P}$ where ${P}$ is a parabolic subgroup of $\mathrm{SO}(4; \C)$. We denote this isomorphism by $\varrho \colon  \mathrm{OG}(1, \C^4) \to \mcal{Q}$ as complex manifolds. 

Regarding $G = \mathrm{SO}(2) = \mathrm{SO}(2;\R) \simeq \mathrm{diag}\left(\mathrm{SO}(2), I_{2} \right)$ as a subgroup of $\mathrm{SO}(4) = \mathrm{SO}(4; \R)$, recall that the quadric $\mcal{Q}$ is acted by $G$ in~\eqref{equ_Gactionqq}. The linear $G$-action on $\C^4$ induces the $G$-action on $\mathrm{OG}(1, \C^4)$. Note that the map $\varrho \colon  \mathrm{OG}(1, \C^4) \to \mcal{Q}$ is $G$-equivariant.

The group $\mathrm{SO}(4)$ acts on the moduli space of one dimensional isotropic subspaces in $\C^4$ linearly and transitively. Also, the adjoint $\mathrm{SO}(4)$-action on $\mcal{O}$ is transitive. We then have 
\begin{itemize}
\item ${\mathrm{SO}(4)}/{\mathrm{S}(\mathrm{O(2)} \times \mathrm{O(2)})} \to \mathrm{OG}(1,\C^4) $
\item $\mcal{O} \to {\mathrm{SO}(4)}/{\mathrm{S}(\mathrm{O(2)} \times \mathrm{O(2)})}$,
\end{itemize}
by choosing one flag and one element of $\mcal{O}$ respectively. We then have a diffeomorphism $\Upsilon \colon \mcal{O} \to  \mathrm{OG}(1, \C^4)$. The orbit $\mcal{O}$ has the adjoint action of the subgroup $G$ of $\mathrm{SO}(4)$. Regarding $G = \mathrm{SO}(2) \simeq \mathrm{diag}\left(\mathrm{SO}(2), I_{2} \right)$ as a subgroup of $\mathrm{SO}(4)$, recall that the quadric $\mcal{Q}$ is acted by $G$ in~\eqref{equ_Gactionqq}. The linear $G$-action on $\C^4$ induces the $G$-action on $\mathrm{OG}(1, \C^4)$. The map $\Upsilon$ is $G$-equivariant. Recall that $\omega_{\mathrm{KKS}}$ is a K\"{a}hler form with respect to the complex structure from $\mathrm{SO}(4; \C)/{P}$.

The orbit $\mcal{O}$ carries two K\"{a}hler forms $(\varrho \circ \Upsilon)^* \omega_{\mcal{Q}}$ and $\omega_{\mathrm{KKS}}$. Each form $\omega^\prime_{t} := (1-t)(\varrho \circ \Upsilon)^* \omega_{\mcal{Q}} + t \omega_{\mathrm{KKS}}$ in the linear interpolation is also symplectic. Since the $G$-action is Hamiltonian with respect to both $\omega^\prime_{0}$ and $\omega^\prime_{1}$, there exists an $G$-equivariant isotopy $\phi^\prime_{t} \colon \mcal{O} \to \mcal{O}$ such that $(\phi^\prime_{t})^* (\omega^\prime_{t}) = \omega^\prime_{0}$ again by Theorem~\ref{thm_equimoserthm}. Set 
$$
\Phi^\prime := \phi^\prime_{1} \circ \Upsilon^{-1} \circ \varrho^{-1} \colon (\mcal{Q},  \omega_{\mcal{Q}}) \to (\mathrm{OG}(1, \C^4), \varrho^* \omega_{\mcal{Q}}) \to (\mcal{O}, (\varrho \circ \Upsilon)^* \omega_{\mcal{Q}}) \to (\mcal{O}, \omega_{\mathrm{KKS}}).
$$
In summary, we have derived the following lemma.

\begin{lemma}\label{lemmaxi2equ}
Under the identification~\eqref{equation_SO2S1}, the map $\Phi^\prime$ is an $S^1$-equivariant symplectomorphism.
\end{lemma}
 
We are now ready to transport the Gelfand--Zeitlin torus $T_{\mathrm{GZ}}$ into $\CP^1 \times \CP^1$. By Lemma~\ref{lemmaxi1equ} and Lemma~\ref{lemmaxi2equ}, we obtain an $S^1$-invariant Lagrangian torus $(\Phi^\prime \circ \Phi)^{-1} (T_{\mathrm{GZ}})$ in the product space $(\CP^1 \times \CP^1, 2(\omega_{\mathrm{FS}} \oplus \omega_{\mathrm{FS}}))$. Let 
\begin{equation}\label{equ_TprimeGZ}
T^\prime_{\mathrm{GZ}}:= (\Phi^\prime \circ \Phi)^{-1} (T_{\mathrm{GZ}}).
\end{equation}

The next proposition compares $T^\prime_{\mathrm{GZ}}$ with $T_{\mathrm{CS}}$ in the same space $(\CP^1 \times \CP^1, 2(\omega_{\mathrm{FS}} \oplus \omega_{\mathrm{FS}}))$. 

\begin{proposition}\label{proposition_gztcstham}
Two monotone Lagrangian tori $T^\prime_{\mathrm{GZ}}$ and $T_{\mathrm{CS}}$ are Hamiltonian isotopic.
\end{proposition}

Let us do preliminary work for verifying Proposition~\ref{proposition_gztcstham}.

Let $\mcal{U}$ be an affine chart of $\CP^1 \times \CP^1$ defined as 
\begin{equation}\label{equ_affinechartU}
\mcal{U}:= (\CP^1 \times \CP^1) \backslash \{v_0 w_0 = 0\} = \CP^1 \backslash (\{v_0 = 0\}) \times  \CP^1 \backslash (\{w_0 = 0\}) \simeq \C \times \C.
\end{equation}
Recall that $\mcal{U}$ is exactly the image of $\rho$ in~\eqref{equ_symplecticembeddingrho}. The Chekanov--Schlenk monotone torus $T_{\mathrm{CS}}$ is contained in the level set $\mu^{-1}_{S^1}(0)$ where the moment map $\mu_{S^1}$ is in~\eqref{equ_momentmapz}. Note that the $S^1$-action generated by $\mu_{S^1}$ induces that on $\mcal{U} \cap \mu^{-1}_{S^1}(0)$.

\begin{figure}[h]
	\scalebox{1.0}{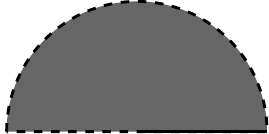}
	\caption{\label{fig_fudamdo} The fundamental domain $\mcal{F}(\sqrt{2})$ in~\eqref{equ_fundamentaldomain}}	
\end{figure}

A fundamental domain of the $S^1$-action on $\mcal{U} \cap \mu^{-1}_{S^1}(0)$ can be chosen as the image of
\begin{equation}\label{equ_fundamentaldomain}
\mcal{F}(\sqrt{2}) := \{\zeta \in \C \mid \mathrm{Im}(\zeta) > 0, \, 0 < |\zeta| < \sqrt{2} \} \cup \{ \zeta \in \C \mid \mathrm{Im}(\zeta) = 0, \, 0 \leq \zeta < \sqrt{2} \}
\end{equation}
under the composition $\rho \circ \Delta_{\mathbb{H}}$ of~\eqref{equ_diagmaph} and~\eqref{equ_symplecticembeddingrho} as depicted in Figure~\ref{fig_fudamdo}. Since each pair of the antipodal points $\zeta$ and $-\zeta$ in $\mathbb{D}(\sqrt{2}):= \{\zeta \in \C \mid 0 \leq |\zeta| < \sqrt{2} \}$ lies on the same $S^1$-orbit, $\mcal{F}(\sqrt{2})$ can be regarded as $\mathbb{D}(\sqrt{2})$ via the map
$$
\Theta \colon \mcal{F}(\sqrt{2}) \to \mathbb{D}(\sqrt{2}), \quad \left( re^{\sqrt{-1} \theta} \mapsto re^{2 \sqrt{-1} \theta} \right).
$$

\begin{lemma}\label{lemma_zerolevelset}
Suppose that a symplectic manifold $(X, \omega)$ admits a Hamiltonian free $S^1$-action. Choose a moment map $\mu_{S^1}$ of the $S^1$-action. Let $L$ be an $S^1$-invariant connected Lagrangian submanifold. Then the Lagrangian submanifold $L$ is contained in the set $\mu^{-1}_{S^1}(r)$ for some $r \in \R \simeq \frak{t}^*$.
\end{lemma}

\begin{proof}
For each point $x \in X$, we know that 
$$
\mathrm{ker}\left(d\mu_{S^1, x}\right) = \left\{ v \in T_x X \mid \omega_x (v, w) = 0 \mbox{ for every $w \in T_x \scr{O}_x$} \right\}
$$
where $\scr{O}_x$ is the orbit through $x$, see \cite[Lemma 5.2.5]{McDuffSalamonIntro}. Since $S^1$ acts freely on the level set $\mu_{S^1}^{-1}(r)$, each point $x \in \mu_{S^1}^{-1}(r)$ is regular and hence $r \in \R$ is a regular value of $\mu_{S^1}$. The dimension counting yields that $\mathrm{ker}\left(d\mu_{S^1, x}\right) = T_x \left(\mu^{-1}_{S^1}(r)\right)$. Since $L$ is Lagrangian, we have
\begin{equation}\label{equ_tangentdescrip}
T_x L \subset  \left\{ v \in T_x X \mid \omega_x (v, w) = 0 \mbox{ for every $w \in T_x \scr{O}_x$} \right\} =  T_x \left(\mu^{-1}_{S^1}(r)\right)
\end{equation}
where $r := \mu_{S^1}(x)$.

We claim that the image of a connected submanifold $L$ under $\mu_{S^1}$ is a singleton set. Suppose on the contrary that the image of $L$ contains an interval having positive measure. By Sard's Theorem, there exist $x \in L$ and $v \in T_x L$ such that $d\mu_{S^1,x}(v)$ is non-zero. Then a non-zero vector $w \in T_x \scr{O}_x$ satisfies $\omega_x (v, w) \neq 0$, which contradicts to~\eqref{equ_tangentdescrip}. 
\end{proof}

\begin{corollary}\label{cor_zerolevelset}
The monotone $S^1$-invariant Lagrangian torus $T^\prime_{\mathrm{GZ}}$ in~\eqref{equ_TprimeGZ} is contained in the zero level set $\mu^{-1}_{S^1}(0)$ where $\mu_{S^1}$ is~\eqref{equ_momentmapz}. 
\end{corollary}

\begin{proof}
Observe that $S^1$ acts on the complement of the four points 
$$
([1:0],[1:0]), ([0:1],[1:0]), ([1:0],[0:1]), \mbox{ and } ([0:1],[0:1])
$$ 
in $\CP^1 \times \CP^1$ freely. Applying Lemma~\ref{lemma_zerolevelset} to the complement, the torus $T^\prime_{\mathrm{GZ}}$ is contained in $\mu^{-1}_{S^1}(r)$ for some $r \in \R$. 

We claim that $T^\prime_{\mathrm{GZ}}$  is contained in the level set of $0$, that is, $r = 0$. By Lemma~\ref{lemma_GZmonotone}, $T_{\mathrm{GZ}}$ is monotone and so is $T^\prime_{\mathrm{GZ}}$ in $\CP^1 \times \CP^1$. Take a point in $T^\prime_{\mathrm{GZ}}$ and choose an $S^1$-invariant almost complex structure $J$. Let $\gamma$ be the integral curve starting from the chosen point of the gradient vector field generated by the Riemannian metric obtained by $\omega$ and $J$. The curve $\gamma$ converges to a certain fixed point. The $S^1$-orbit of $\gamma$ is a $J$-holomorphic disk such that its symplectic area is equal to $\pi / 2$ times its Maslov index. Such a $J$-holomorphic disk is bounded by an $S^1$-orbit in $T_{\mathrm{GZ}}$ and called a \emph{gradient holomorphic disk}, see \cite[Section 2]{ChoKimMONO} for more details. As the level of $T^\prime_{\mathrm{GZ}}$ varies, the symplectic area of this gradient disk changes due to Archimedes, while the Maslov index of the disk does not change. It implies that the monotone torus must be contained in the level zero. 
\end{proof}

Corollary~\ref{cor_zerolevelset} says that the torus $T^\prime_{\mathrm{GZ}}$ is completely contained in the $S^1$-orbit of the image of the open disk $\mathbb{D}(\sqrt{2}) $ under the composition map $\rho \circ \Delta_{\mathbb{D}}$ in~\eqref{equ_diagmap} and~\eqref{equ_symplecticembeddingrho}. By intersecting the fundamental domain, we obtain a curve $\Gamma_{\mathrm{GZ}}$ determined by 
$$
\Gamma_{\mathrm{GZ}} := \Theta \left( (\rho \circ \Delta_\mathbb{D})^{-1} (T^\prime_{GZ}) \cap \mcal{F}(\sqrt{2}) \right).
$$
Since $T^\prime_{GZ}$ is a smooth torus,  $\Gamma_{\mathrm{GZ}}$ is a simple closed curve not passing through the origin. Then there are two possibilities$\colon$
\begin{enumerate}
\item $\Gamma_{\mathrm{GZ}}$ does not bound the closed region,
\item $\Gamma_{\mathrm{GZ}}$ bounds the closed region.
\end{enumerate}
We claim that $(2)$ holds.

\begin{proposition}\label{lemma_gammagzenlose}
The curve $\Gamma_{\mathrm{GZ}} \subset \mathbb{D}(\sqrt{2})\backslash \{0\}$ bounds the closed region.
\end{proposition}

To see Proposition~\ref{lemma_gammagzenlose}, we collect some facts. Let $(X, \omega)$ be a closed monotone symplectic manifold and $L$ a monotone Lagrangian submanifold. For each homotopy class $\beta \in \pi_2(X,L)$ of Maslov index two, a generic choice of $\omega$-compatible almost complex structure makes the moduli space of stable maps in the class $\beta$ from $(\mathbb{D}, \partial \mathbb{D})$ to $(X,L)$ transversal. The virtual dimension of this moduli space is exactly dimension of $L$. This monotonicity condition ensures that the moduli space is closed. We then count the number of stable map passing through a generic point of $L$ at the marking point on the boundary of the disk $\mathbb{D}$. Such a number is called a \emph{counting invariant} of $L$ in $X$ and is denoted by $n_\beta$. The sum of counting invariants is meant to be $\sum_\beta n_\beta$ where the summation is taken over all homotopy class $\beta \in \pi_2(X,L)$. By the dimension reason, almost all $n_\beta = 0$ and the sum is finite because of Gromov's compactness theorem.

\begin{lemma}[\cite{EliashbergPolterovich}]\label{lemma_countings}
Suppose that $\phi$ is a symplectomorphism from $(X, \omega)$ to $(X^\prime, \omega^\prime)$. If $L$ and $L^\prime$ are monotone Lagrangian submanifolds related by $\phi$, then the sum of counting invariants bounded by $L$ in $X$  is equal to that of counting invariants by $L^\prime$ in $X^\prime$.
\end{lemma}

The sum of counting invariants of $T_{\mathrm{GZ}}$ in $\mcal{O}$ can be computed from the disk potential of $T_{\mathrm{GZ}}$ in \cite[Theorem A]{KimQuad}.

\begin{lemma}[\cite{KimQuad}]\label{lemma_sumofcounting}
The sum of counting invariants of $T_{\mathrm{GZ}}$ in $\mcal{O}$ is five.
\end{lemma} 

\begin{proof}[Proof of Proposition~\ref{lemma_gammagzenlose}]
Suppose that the curve $\Gamma_{\mathrm{GZ}}$  does not bound the closed region. Then the torus $T^\prime_{\mathrm{GZ}}$ is Hamiltonian isotopic to the product of equators of $S^2 \times S^2$. If so, the sum of counting invariants of $T^\prime_{\mathrm{GZ}}$ in $S^2 \times S^2$ is four. According to Lemma~\ref{lemma_countings} and Lemma~\ref{lemma_sumofcounting},  the sum of counting invariants of $T^\prime_{\mathrm{GZ}}$ has to be five. We have derived a contradiction.
\end{proof}

Recall that the chosen $\Gamma$ for~\eqref{equ_TCS} is contained in $\mathbb{H}(\sqrt{2})$. Let 
$$
\Gamma_{\mathrm{CS}} := \Theta(\Gamma).
$$

Let us compare two simple closed curves $\Gamma_{\mathrm{CS}}$ and $\Gamma_{\mathrm{GZ}}$. A priori, the simple closed curve $\Gamma_{\mathrm{GZ}}$ might be complicated so that it may not be contained in a branch, while $\Gamma_{\mathrm{CS}}$ is contained in the branch $\mathbb{D}(\sqrt{2}) \backslash \{ z \mid \mathrm{Im}(z) = 0, z \geq 0\}$. 

We shall deform one to the other via a Hamiltonian isotopy. For this purpose, we recall one lemma concerning Hamiltonian isotopy class of loops in a two dimensional exact symplectic manifold. (Here, the dimension condition is necessary to apply the Moser argument).

\begin{lemma}[Lemma 2.3 in \cite{LekiliMaydanskiy}]\label{lemma_aa}
Let $(\Sigma, \omega := d \alpha)$ be an exact symplectic manifold of real dimension two. Let $\Gamma_1$ and $\Gamma_2$ be simple closed curves in $\Sigma$.
Suppose that
\begin{enumerate}
\item $\Gamma_0$ and $\Gamma_1$ are isotopic and
\item $\int_{\partial \Gamma_0} \alpha = \int_{\partial \Gamma_1} \alpha$.
\end{enumerate}
Then $\Gamma_0$ and $\Gamma_1$ are related by via a compactly supported Hamiltonian isotopy. 
\end{lemma}

In our case, two simple closed curves $\Gamma_{\mathrm{CS}}$ and $\Gamma_{\mathrm{GZ}}$ in $\mathbb{D}(\sqrt{2}) \backslash \{0\}$ bound regions having the same area since the bounded regions can be lifted to a disk bounded by $T_{\mathrm{CS}}$ and by $T^\prime_{\mathrm{GZ}}$ respectively. The monotonicity of $T_{\mathrm{CS}}$ and $T^\prime_{\mathrm{GZ}}$ ensures the bounded regions are same. By Stokes' theorem, the second condition holds. For the first condition, we can apply the following well-known fact from differential topology.

\begin{lemma}[Theorem 3.1 in Chapter 8 of \cite{HirschDT}] 
Let $X$ be a connected $n$-manifold and $f,g \colon \mathbb{D}^k \to X$ embedding of the $k$-disk, $0 \leq k \leq n$. If $k = n$ and $X$ is orientable, assume that $f$ and $g$ both preserve orientation. Then $f$ and $g$ are isotopic. 

If $f(\mathbb{D}^k) \cup g(\mathbb{D}^k)$ is contained in $X - \partial X$, an isotopy between them can be realized by a smooth isotopy of $X$ having compact support.
\end{lemma}

We now discuss the relation between $T_{\mathrm{CS}}$ and by $T^\prime_{\mathrm{GZ}}$. A choice of isotopy from $\Gamma_{\mathrm{CS}}$ to $\Gamma_{\mathrm{GZ}}$ in $\mathbb{D}(\sqrt{2}) \backslash \{0\}$ gives rise to a Lagrangian isotopy between $T_{\mathrm{CS}}$ and by $T^\prime_{\mathrm{GZ}}$. Proposition~\ref{proposition_gztcstham} further claim that the Lagrangian isotopy arising from the Hamiltonian isotopy in $\mathbb{D}(\sqrt{2}) \backslash \{0\}$ can be extended to an ambient Hamiltonian isotopy of $\CP^1 \times \CP^1$. 

\begin{proof}[Proof of Proposition~\ref{proposition_gztcstham}]
For simplicity of notation, let us set
$$
\left( X := \CP^1 \times \CP^1, \omega := 2 (\omega_{\mathrm{FS}} \oplus \omega_{\mathrm{FS}}) \right),
$$
while presenting the proof. Suppose that $\ell \colon T^2 \times [0,1]  \to \mathbb{D}(\sqrt{2}) \times \mathbb{D}(\sqrt{2})$ is a Lagrangian isotopy arising from the Hamiltonian isotopy from Lemma~\ref{lemma_aa}. Since $\ell$ is a Lagrangian isotopy, the pull-backed form is of the following form
$$
(\Delta_{\mathbb{D}} \circ \ell)^* \omega = \ell^* \omega_\mathrm{std}=\alpha_s \wedge ds
$$
where $\{\alpha_s\}$ is a family of one forms. 

We shall show that the Lagrangian isotopy $\Delta_\mathbb{D} \circ \ell$ is exact. Then this Lagrangian isotopy can be extended to an ambient Hamiltonian isotopy of $\CP^1 \times \CP^1$ as desired, see \cite[Section 6]{Polterovichgroupof}, \cite[Section 3.6]{OhBook1} for instance. 

It remain to show that for each $s$, $\alpha_s$ is exact. Note that $\omega_{\mathrm{std}} = d \eta$. The pull-back of the primitive $\eta$ is 
$$
\ell^* \eta = f(x,s) ds + \eta_s^\prime
$$
where $\{\eta^\prime_s\}$ is a family of one forms on $T^2 \times \{s\}$.
To show that $\alpha_s$ is exact, we claim that the integration of $\eta$ over any loop remains constant through the isotopy. Let $\Gamma_s$ be an isotopy from $\Gamma_{\mathrm{CS}}$ to $\Gamma_{\mathrm{GZ}}$ from Lemma~\ref{lemma_aa}. We then have an isotopy of loops 
$$
\gamma_s:= (\rho \circ \Delta_\mathbb{D}) \left(\Theta^{-1}(\Gamma_s)\right).
$$ 
Then $\int_{\gamma_s} \eta$ is constant through the isotopy because of the construction of the isotopy.
We denote by $\gamma_0^\prime$ a Lefschetz thimble bounded by $T_{\mathrm{CS}}$ in $\mathbb{D}(\sqrt{2}) \times \mathbb{D}(\sqrt{2})$. Let $\gamma_s^\prime$ be the isotoped circle. In other words,
$\gamma_s^\prime := \left( \ell_s \circ \ell^{-1}_0 \right) \left( \gamma_0^\prime \right)$ where $\ell_s := \ell (\cdot, s) \colon T^2 \to \mathbb{D}(\sqrt{2}) \times \mathbb{D}(\sqrt{2})$. Then $\int_{\gamma^\prime_s} \eta = 0$ for all $s$. Hence, the claim is verified.

The claim yields that the integration of $\eta^\prime_s$ over each loop in $T^2$ is independent to $s$. Then $d_x f  = \alpha_s$ as desired.
\end{proof}

\begin{proof}[Proof of Theorem~\ref{theorem_main}]
It follows from Lemma~\ref{lemmaxi1equ}, Lemma~\ref{lemmaxi2equ}, and Proposition~\ref{proposition_gztcstham}.
\end{proof}

\providecommand{\bysame}{\leavevmode\hbox to3em{\hrulefill}\thinspace}
\providecommand{\MR}{\relax\ifhmode\unskip\space\fi MR }
\providecommand{\MRhref}[2]{%
  \href{http://www.ams.org/mathscinet-getitem?mr=#1}{#2}
}
\providecommand{\href}[2]{#2}

\end{document}